\documentclass{amsart}

\usepackage{amssymb,amscd}
\usepackage{enumerate,comment}
\theoremstyle{plain}

\newtheorem{thm}{Theorem}[section]
\newtheorem{lem}[thm]{Lemma}
\newtheorem{prop}[thm]{Proposition}
\newtheorem{cor}[thm]{Corollary}

\theoremstyle{definition}
\newtheorem{defn}[thm]{Definition}

\theoremstyle{remark}

\def\nul{\emptyset }

\def\d{\delta }

\def\e{\varepsilon }

\def\S{\Sigma }

\newcommand{\mbf}{\mathbf}
\newcommand{\la}{\langle}
\newcommand{\ra}{\rangle}
\newcommand{\NN}{\mathbb{N}}

\newcommand{\cF}{\mathcal{F}}

\newcommand{\cL}{\mathcal{L}}
\newcommand{\cM}{\mathcal{M}}
\newcommand{\cN}{\mathcal{N}}
\newcommand{\cS}{\mathcal{S}}
\newcommand{\Lp}{\mathcal{L}^{\textrm{pre}}}
\newcommand{\Gp}{\mathcal{L}_{\cS}^{\textrm{pre}}}
\newcommand{\Kp}{\mathcal{L}_{K}^{\textrm{pre}}}

\newcommand{\Th}{\texttt{Th}}
\newcommand{\The}{\texttt{Th}_\exists}

\newcommand{\factor}[2]{{\raise0.7ex\hbox{$#1$} \!\mathord{\left/
 {\vphantom {#1 {#2}}}\right.\kern-\nulldelimiterspace}
\!\lower0.7ex\hbox{${#2}$}}}

\newcommand{\be}{\begin{enumerate}}
\newcommand{\ee}{\end{enumerate}}

\newcommand{\mcomment}[1]{\marginpar{\footnotesize #1}}

\title{Stability of Universal Equivalence of Groups under Free Constructions}

\author[A. Duncan]{Andrew J. Duncan}
\address{School of Mathematics and Statistics, University of Newcastle-upon-Tyne, Newcastle-Upon-Tyne NE1 7RU, United Kingdom}
\email{a.duncan@ncl.ac.uk}

\author[I. Kazachkov]{Ilya V. Kazachkov}
\address{Department of Mathematics and Statistics, McGill University,
805 Sherbrooke St. West, Montreal, Quebec H3A 2K6, Canada}
\email{kazachkov@math.mcgill.ca}

\author[V. Remeslennikov]{Vladimir N. Remeslennikov}
\address{Institute of Mathematics (Russian Academy of Science),
13 Pevtsova St., Omsk, 644099, Russia}
\email{remesl@iitam.omsk.net.ru}

\thanks{Research supported by EPSRC grant EP/D065275}

\begin{document}
\maketitle

\section{Introduction}

In his important paper in \cite{S} J. Stallings introduced a
generalisation of amalgamated products of groups -- called a pregroup, which is a particular kind of a partial group. He then defined the universal group $U(P)$ of a pregroup $P$ to be a universal object (in the sense of category theory) extending the partial operations on $P$ to group operations on $U(P)$. The universal group turned out to be a versatile and convenient generalisation of classical group constructions: HNN-extensions and amalgamated products. In this respect the following general question arises.Which properties of pregroups, or relations between pregroups,  carry over to the respective universal groups? The aim of this paper is to prove that universal equivalence of pregroups extends to universal equivalence of their universal groups.

We begin by some preliminary model-theory results. We refer the reader to \cite{marker} for a detailed introduction to model theory. The main goal here is to give a criterion of universal equivalence of two models in the form that best suits our needs.

\section{Preliminaries} \label{sec:pre}

Let $\cL$ be a language with signature $(C, F, R)$, and variables $X$,
where $C$ is a set of constants and $F$ and $R$ are  finite sets of
functions and relations respectively. In addition each element
$f$ of $F$ is
associated to a non-negative integer  $n_f$, and similarly for $R$.

An $\cL$-structure $\cM$ is a 4-tuple:
\begin{itemize}
  \item a non-empty set $M$;
  \item a function $f_\cM: M^{n_f} \to M$ for each $f\in F$;
  \item a set $r_\cM\subseteq M^{n_r}$ for each $r\in R$;
  \item an element $c_\cM$ for each element $c\in C$.
\end{itemize}
If $a$ is an element of $C$, $F$ or $R$ we refer to
$a_\cM$ as the interpretation of $a$
in $\cM$.
The subscript $\cM$ is omitted where no ambiguity arises.

We use the language $\cL$ to write formulas describing the properties of $\cL$-structures. Roughly speaking formulas are constructed inductively starting from constant symbols from $C$ and variable symbols $v_1,\dots, v_n,\dots$, using the Boolean connectives, relations from $R$, functions from $F$ and the equality symbol `$=$'.

More precisely,
the set of $\cL${\em -terms} is the smallest set $T$ such that:
\begin{itemize}
  \item $c\in T$ for each constant symbol $c\in C$;
  \item each variable symbol $v_i\in T$;
  \item if $t_1,\dots, t_n\in T$ and $f\in F$
then $f(t_1,\dots, t_{n_f})\in T$.
\end{itemize}

We say that $\Phi$ is an {\em atomic} $\cL${\em -formula} if $\Phi$ is either
\begin{itemize}
  \item $t_1=t_2$, where $t_1$ and $t_2$ are terms or,
  \item $r(t_1,\dots, t_{n_r})$, where $r\in R$ and $t_1,\dots, t_{n_r}$
are terms.
\end{itemize}

The set of $\cL${\em -formulas} is the smallest set $W$ containing atomic formulas and such that
\begin{itemize}
\item if $\Phi\in W$ then $\neg \Phi\in W$;
\item if $\Phi$ and $\Psi$ are in $W$ then $\Phi \wedge \Psi$ and $\Phi \vee \Psi$ are in $W$ and
\item if $\Phi$ is in $W$ then  $\exists v_i \Phi$ and $\forall v_i \Phi$ are in $W$.
\end{itemize}

It is often useful in practice to observe that, as
\[\Phi \vee \Psi\equiv \neg (\neg \Phi \wedge \neg \Psi)\]
and
\[\forall v_i \Phi \equiv \neg \exists v_i(\neg \Phi)\]
we can construct all formulas (up to logical equivalence $\equiv$) without using
$\vee$ or $\forall$.

To make induction arguments precise we shall define, for any term or formula $s$ of $\cL$, the {\em level} $l(s)$ and  the {\em constants} $C(s)$ of $s$. If $s$ is a formula we shall also define the  {\em degree} $d(s)$ of $s$. To begin with if $t$ is a term and  $t=x$ or $t=c$, where $x$ is a variable and $c$ a constant, then $l(t)=0$ and
\[
C(t)=
\left\{
\begin{array}{ll}
\nul,& \textrm{ if }t=x\\
c,& \textrm{ if } t=c
\end{array}
\right.
.
\]
If $t=f(t_1,\ldots ,t_{n})$,
where $n=n_f$ and the $t_i$ are terms then
$l(t)=\max\{l(t_1),\ldots, l(t_n)\}+1$,
and $C(t)=
\cup_{i=1}^n C(t_i)$.

If $a$ is an atomic formula of the form $t_1=t_2$, for terms $t_1$ and $t_2$,
then we
define $l(a)=\max\{l(t_1),l(t_2)\}$
and $C(t)=C(t_1)\cup C(t_2)$.
If $r$ is an $n$-ary relation and $a=r(t_1,\ldots, t_n)$ then
set $l(a)= \max\{l(t_1),\ldots ,l(t_n)\}$
and $C(a)=\cup_{i=1}^n C(t_i)$.
The degree of an atomic formula $a$ is defined to be $d(a)=0$.

If $\Phi=\neg \Psi$ or $\Phi =\exists x \Psi$ or $\forall x \Psi$
then we define
$l(\Phi)=l(\Psi)$, $d(\Phi)=d(\Psi)+1$
and $C(\Phi)=C(\Psi)$.
If $\Phi= \Phi_1\wedge \Phi_2$ or $\Phi_1\vee \Phi_2$ then we define
$l(\Phi)=\max\{l(\Phi_1),l(\Phi_2)\}$, $d(\Phi)=d(\Phi_1)+d(\Phi_2)$
and $C(\Phi)=C(\Phi_1)\cup C(\Phi_2)$.

We say that a variable $v$ occurs freely in a formula $\Phi$ if it is not inside a $\exists v$ or a $\forall v$ quantifier, otherwise $v$ is said to be bound. A formula is called a {\em sentence} or {\em closed} 
if it has no free variables.

Let $\Phi$ be a formula with free variables from $v=(v_1,\dots, v_m)$ and let $\bar a=(a_1,\dots, a_m)\in M^m$. We inductively define when $\Phi$
 {\em holds} on $\bar a$ in a $\cL$-structure $\cM$ ($\Phi(\bar a)$ is
{\em true} in $\cM$ or $\cM$ {\em satisfies} $\Phi(\bar a)$), write $\cM\models \Phi(\bar a)$.
\begin{itemize}
  \item if $\Phi$ is $t_1=t_2$, then $\cM \models \Phi (\bar a)$ if $t_1(\bar a)=t_2(\bar a)$;
  \item if $\Phi=r(t_1,\dots,t_{n_r})$, then $\cM \models \Phi(\bar a)$ if $r(t_1(\bar a),\dots,t_{n_r}(\bar a))\in r_\cM$;
  \item if $\Phi=\neg \Psi$ then $\cM \models \Phi(\bar a)$ if $\cM\nvDash \Psi(\bar a)$;
  \item if $\Phi=\Psi_1\wedge \Psi_2$ then $\cM \models \Phi(\bar a)$ if $\cM\models \Psi_1(\bar a)$ and $\cM\models \Psi_2(\bar a)$;
  \item if $\Phi=\Psi_1\vee \Psi_2$ then $\cM \models \Phi(\bar a)$ if $\cM\models \Psi_1(\bar a)$ or $\cM\models \Psi_2(\bar a)$;
  \item if $\Phi =\exists v_{m+1} \Psi(\bar v, v_{m+1})$, then $\cM\models \Phi$ if there exists $b\in M$ such that $\cM \models \Psi(\bar a, b)$;
  \item if $\Phi =\forall v_{m+1} \Psi(\bar v, v_{m+1})$, then $\cM\models \Phi$ if for all $b\in M$ one has $\cM \models \Psi(\bar a, b)$
\end{itemize}

A set of sentences is called a \emph{theory}. We say that $\cM$ is a
{\em model} of a theory $T$ if $\cM\models \Phi$ for all $\Phi\in T$. For an $\cL$-structure $\cM$ we denote by $\Th(\cM)$ the collection of all sentences that are satisfied by $\cM$, $\Th(\cM)$ is called the {\em full} or
{\em elementary theory} of $\cM$.

Every formula $\Phi$ of $\cL$ with free variables $\bar v =(v_1,\ldots ,v_k)$  is logically equivalent to a formula of the type
$$
Q_1x_1 Q_2 x_2 \ldots Q_n x_n \Psi(\bar x,\bar v),
$$
where  $Q_i \in \{\forall, \exists \}$, and  $\Psi(\bar x,\bar v)$ is a boolean combination of atomic formulas in
variables from $\bar v \cup \bar x$. This form is called the
{\em prenex normal form} of a formula $\Phi$.

A sentence $\Phi$ is called {\em universal} ({\em existential}) if $\Phi$ is equivalent to a formula of the form
$$
Q_1x_1 Q_2 x_2 \ldots Q_n x_n \Psi(\bar x),
$$
where  $Q_i = \forall$ ($Q_i=\exists$) for all $i$, and  $\Psi(\bar x)$ is a boolean combination of atomic formulas in the indicated variables.
The collection of all universal (existential) sentences that are satisfied by an $\cL$-structure $\cM$ is called the {\em universal} ({\em existential}) theory of $\cM$, we denote it by $\Th_\forall (\cM)$ ($\Th_\exists (\cM)$).
If $\cM$ and $\cN$ are $\cL$-structures and
$\Th_\forall (\cM)=\Th_\forall (\cN)$ we say that $\cM$ and $\cN$ are
{\em universally equivalent} and write $\cM\equiv_\forall \cN$.
{\em existential equivalence} is defined similarly and we write
$\cM\equiv_\exists \cN$ if $\cM$ and $\cN$ are
existentially equivalent.

Let $A$ be a set of sentences
of $\cL$ and let $\cM$ and $\cN$ be models of $A$
with underlying sets $M$ and $N$ respectively.
For subsets $S$ and $T$
of $M$ and $N$ respectively we say that a map $\phi:S\rightarrow T$ is
an $\cL${\em -morphism} if the following conditions hold.
\be
\item If $c\in C\cap S$ then $c\in T$ and $\phi(c)=c$.
\item If $f$  is an $n$-ary function (i.e. $n_f=n$) in $F$ and $f(s_1,\ldots ,s_n)\in S$,
for some $n$-tuple $(s_1, \ldots ,s_k)$ of elements of $S$, then
\[\phi( f(s_1,\ldots ,s_n) )= f(\phi(s_1),\ldots, \phi(s_n))\in T.\]
\item If $r$  is an $n$-ary relation in $R$ and $(s_1,\ldots ,s_n)\in r$,
for some $n$-tuple $(s_1, \ldots ,s_n)$ of elements of $S$, then
$(\phi(s_1),\ldots, \phi(s_n))\in r$.
\ee
If $\phi:S\rightarrow T$ is a bijective $\cL$-morphism
such that $\phi^{-1}$ is an
$\cL$-morphism from $T$ to $S$ then we say that $\phi$ is an
$\cL${\em -isomorphism} and that $S$ and $T$ are $\cL${\em -isomorphic}
or $S\cong_{\cL} T$.
$\cL$-isomorphism defines an equivalence relation on the subsets of a
model $\cM$ and we denote by $[S]$ the equivalence class of $S$.

Now we restrict attention to finite subsets of models. We denote by
$\cF_\cL(\cM)=\cF(\cM)$ the set of $\cL$-isomorphism equivalence classes of
finite subsets of $M$. We say that models $\cM$ and $\cN$
have {\em equivalent} $\cL$-isomorphism
classes of finite subsets, and write $\cF(\cM)\equiv \cF(\cN)$, if there exists
a bijection $\theta : \cF(\cM)\rightarrow \cF(\cN)$ such that, for all
finite subsets $S\subseteq M$, if $\theta ([S]) = [T]$ then
there exists an $\cL$-isomorphism $\phi(S)\rightarrow T^\prime$, for
some $T^\prime \in [T]$ (hence for all $T^\prime \in [T]$).
\begin{lem}\label{lem:pseq}
$\cF(\cM)\equiv \cF(\cN)$ if and only if, for all finite subsets $S\subseteq M$,
there exists a subset $T\subseteq N$ such that $S\cong_\cL T$.
\end{lem}
\begin{proof}
If  $\cF(\cM)\equiv \cF(\cN)$ and $S$ is a finite subset of $M$ then,
by definition, $S$ is $\cL$-isomorphic to some finite subset of $N$.
Conversely, suppose every finite subset of $M$ is $\cL$-isomorphic
to a finite subset of $N$.  For each  isomorphism class $U$ of finite
subsets of $M$ choose a representative $S_U$, so $U=[S_U]$. Similarly choose
a representative $T_V$ for each isomorphism class of finite subsets
of $N$. Consider an isomorphism class $U\in \cF(\cM)$. $S_U$ is
$\cL$-isomorphic  to $T$ for some finite subset of $N$. Let $T^\prime$ be
the chosen representative of $[T]$. Then $S_U\cong_\cL T^\prime$. Define
$\theta (U)=[T^\prime]$. Then $\theta$ is a well-defined map from
$\cF(\cM)$ to $\cF(\cN)$ and straightforward verification shows that
$\theta $ is a bijection. By construction, if $\theta(U)=V$ then $S_U$
is $\cL$-isomorphic to the representative $T^\prime$ of $V$,
so the same goes of any element $S\in U$.
Hence $\cF(\cM)\equiv \cF(\cN)$.
\end{proof}
If $\cM$ and $\cN$ are models of $A$ then it is easy to see that
$\cM\equiv_\exists \cN$ if and only if $\cM\equiv_\forall \cN$. The
following proposition gives a further characterisation of this property,
in certain cases.

\begin{prop}\label{prop:fe}
Assume that $\cL$ has signature $(C,F,R)$ where either
\be
\item $C$ is finite or
\item
$R$ contains a relation $\d_C$ and, for each $c\in C$, $A$ contains
axioms
\be
\item
$c\in \d_C$ and
\item
$\forall x (x\notin \d_C \implies x\neq c)$.
\ee
\ee
Let $\cM$ and $\cN$ be models of $A$. Then $\cM\equiv_\exists \cN$ if and
only if $\cF(\cM)\equiv \cF(\cN)$.
\end{prop}
\begin{proof}
Assume first that $\cM\equiv_\exists \cN$. Let $S=\{m_1,\ldots, m_k\}$ be
 a finite subset of $\cM$. Define the formula
\[\Phi_1=\bigwedge_{1\le i<j\le k} x_i\neq x_j.\]
$\Phi_1$ will enable us to identify  $k$ distinct elements of
$M$ or $N$; in fact $\cM\vDash \Phi_1[m_1,\ldots ,m_k]$ so $\cM\vDash \exists x_1,\ldots ,x_k
\Phi_1$.

Now let $S\cap C=\{m_{i_1},\ldots ,m_{i_s}\}$, say $m_{i_j}=c_j\in C$
and let $\{1,\ldots, k\}\backslash \{i_1,\ldots, i_s\}=\{j_1,\ldots, j_t\}$.
Define
\[\Phi_2=\left(\bigwedge_{r=1}^s x_{i_r}=c_r \right)\wedge
\left(\bigwedge_{c\in C} \bigwedge_{r=1}^t x_{j_r}\neq c\right),\]
if $C$ is finite and
\[\Phi_2=\left(\bigwedge_{r=1}^s x_{i_r}=c_r \right)\wedge
\left(\bigwedge_{r=1}^t x_{j_r}\notin \d_C\right),\]
otherwise. By construction $\cM\vDash \Phi_2[m_1,\ldots ,m_k]$ and
$\Phi_2$ allows us to identify $C\cap S$ and a corresponding
subset of $N$.

Write $I_k=\{1,\ldots ,k\}$.
Let $f\in F$ be an $n$-ary function, for some $n\ge 1$. Let
\[S_{f,0}=
\{(i_1,\ldots, i_n)\in I_k^n| f(m_{i_1},\ldots ,m_{i_n})\in S\backslash C\}\]
and
\[S_{f,1}=
\{(i_1,\ldots, i_n)\in I_k^n| f(m_{i_1},\ldots ,m_{i_n})\in C\}.\]
For each $(i_1,\ldots, i_n)\in S_{f,1}$ define $s=s(i_1,\ldots ,i_n)$
to be the integer in $I_k$ such that $f(m_{i_1},\ldots ,m_{i_n})=
m_s$.
Define
\[\Phi_{f,0} =\bigwedge_{(i_1,\ldots, i_n)\in S_{f,0}} f(x_{i_1},\ldots ,x_{i_n})=f(m_{i_1},\ldots ,m_{i_n})\]
and
\[\Phi_{f,1} =
\bigwedge_{(i_1,\ldots, i_n)\in S_{f,1}}
f(x_{i_1},\ldots ,x_{i_n})=m_{s(i_1,\ldots ,i_n)}
.
\]
Define $\Phi_f=\Phi_{f,0}\wedge \Phi_{f,1}$. Then $\cM\vDash
\Phi_f[m_1,\ldots ,m_k]$.

Let $r\in R$ be an $n$-ary relation, for some $n\ge 1$,
 and let
\[S_r=\{(i_1,\ldots ,i_n)\in I_k| (m_{i_1},\ldots ,m_{i_n})\in r\}.\]
Define
\[
\Phi_r = \left(\bigwedge_{(i_1,\ldots ,i_n)\in S_r }
r(x_{i_1},\ldots, x_{i_n})
 \right)
\wedge
\left(\bigwedge_{(i_1,\ldots ,i_n)\notin S_r }
\neg r(x_{i_1},\ldots, x_{i_n})
 \right)
\]
Then $\cM\vDash
\Phi_r[m_1,\ldots ,m_k]$.

Finally define $\Phi=\Phi_1\wedge \Phi_2 \wedge \bigwedge_{f\in F}
\Phi_f \wedge \bigwedge_{r\in R}\Phi_r$. Then
$\cM\vDash \Phi[m_1,\ldots, m_k]$ so $\cM\vDash \exists x_1,\ldots, x_k
\Phi$. Therefore $\cN\vDash \exists x_1,\ldots, x_k
\Phi$ and there exist $n_1,\ldots ,n_k \in N$ such that
$\cN\vDash \Phi[n_1,\ldots, n_k]$.

Set $T=\{n_1,\ldots ,n_k\}$
and define $\phi:S\rightarrow T$ by $\phi(m_i)=n_i$, $i=1,\ldots ,k$.
By definition $\phi$ is an $\cL$-morphism and is a bijection of $S$ and
$T$. Moreover $\phi^{-1}$ is, by construction of $\Phi$, an $\cL$-morphism.
Hence $S\cong_\cL T$ and it follows from Lemma \ref{lem:pseq} that $\cF(\cM)
\equiv \cF(\cN)$.

Now suppose that $\cF(\cM)
\equiv \cF(\cN)$.
Write $F_n$ for the set of $n$-ary functions of $F$.
Since $F$ is finite we may assume that $F$ is the union of
$F_n$, for $n$ from $1$ to
$K$, for some $K\in\NN$.
 Given a finite subset $S$ of $M$ we define the following sequence
of subsets.
Set $S_0=S$ and having defined $S_i$ set
\[S_{i+1}=S_i\cup \bigcup_{n=1,\ldots ,K}\bigcup_{f\in F_n}
\{f(m_1,\ldots, m_n)|m_j\in S_i, j=1,\ldots ,n\}.
\]
Now choose $T_l\subseteq N$ such that there is an $\cL$-isomorphism
$\phi_l$ from $S_l$ to $T_l$, for all $l\ge 0$.

Consider a term $t$ of level $l$ with variables among $x_1,\ldots, x_k$
and a $k$-tuple $a_1,\ldots, a_k$ of elements of $M$ and set
$S=\{a_1,\ldots ,a_k\}\cup C(t)$. We claim that $t(a_1,\ldots, a_k)\in
S_l$. To see this note that it holds when $l=0$, since in this
case $t(a_1,\ldots, a_k)\in S$. Suppose then that $t$ has level
$l$ and that the claim holds for at all levels below $l$. Then
$t=f(t_1,\ldots, t_m)$, where $l(t_i)<l$. By assumption
$t_i(a_1,\ldots, t_k)\in S_{l-1}$ and so by definition $t(a_1,\ldots, a_k)
=f(t_1(a_1,\ldots, t_k),\ldots ,t_m(a_1,\ldots, t_k))\in S_l$;
and the claim holds for all $l$ by induction.

Let $\Phi$ be a quantifier free formula with variables among
$x_1,\ldots, x_m$ and let
$\Psi=\exists x_1,\ldots ,x_m \Phi$. We wish to show that
$\cM\vDash \Psi$ if and only if $\cN\vDash
\Psi$. To do this we shall proceed as follows. Suppose
 $\Phi$ has level $l$ and let $a_1,\ldots ,a_m\in M$.
Let  $S=S(\Phi)=\{a_1,\ldots, a_m\}\cup C(\Phi)$, where $C(\Phi)$ is
the set of constants of $\Phi$, and define $S_0,S_1,\ldots$ and
$\phi_1,\phi_2,\ldots$ as above. We shall prove that
\begin{equation}\label{eq:indhyp}
\cM\vDash \Phi(a_1,\ldots ,a_m)
\textrm{
if and only if
 }
\cN\vDash \Phi(\phi_l(a_1),\ldots ,\phi_l(a_m))
,
\end{equation}
and the
result will follow immediately.
We use induction on $(d,l)$,
where $d$ is the degree and $l$ the level of $\Phi$.

Assume that
\eqref{eq:indhyp} holds whenever $\Phi$ has level at most $l$ and degree
$0$. Suppose now that
$\Phi$ has level $l+1$ and  degree $0$.
 In this case $\Phi$ is of the form $t_1=t_2$,
or of the form $r(t_1,\ldots, t_m)$,
where $r\in R$ and the $t_i$ are terms. Since $\Phi$ has level $l+1$ at
least one of the $t_i$ has level $l+1$ and none have level greater
than $l+1$. Hence $t_i(a_1,\ldots, a_k)\in S_{l+1}$, for all $i$.
Set $b_i=\phi_{l+1}(a_i)$, $i=1,\ldots ,k$.
Then, as $\phi_{l+1}$ is an isomorphism with domain $S_{l+1}$,
we have $\phi_{l+1}(t_i(a_1,\ldots ,a_k))
=t_i(b_1,\ldots ,b_k)$, for all $i$. Furthermore $t_1(a_1,\ldots, a_k)=
t_2(a_1,\ldots, a_k)$ if and only if
$t_1(b_1,\ldots ,b_k)=t_2(b_1,\ldots, b_k)$
and $(t_1(a_1,\ldots, a_k),\ldots ,t_m(a_1,\ldots, a_k))\in r$ if and only
if  $(t_1(b_1,\ldots, b_k),\ldots ,t_m(b_1,\ldots, b_k))\in r$. Hence
$\cM\vDash \Phi(a_1,\ldots a_k)$
if and only if $\cN\vDash \Phi(b_1,\ldots,b_k)$. Therefore the result
holds for $\Phi$ of level $l+1$ and degree $0$. Note that this argument
also goes through in the case $(d,l)=(0,0)$ so by induction
\eqref{eq:indhyp} holds for formulae $\Phi$ of level $l$ and degree
$0$, for all non-negative integers $l$.

Now let $d$ and $l$ be non-negative integers and assume that
\eqref{eq:indhyp} holds for formulae $\Phi$ of degree $d_1$ and
level $l_1$ where  either (i) $d_1\le d$
and  $l_1=l$ or (ii) $l_1<l$.
Suppose then that  $\Phi$ has level $l$ and degree $d+1$.
Then
either $\Phi =\neg \Phi_1$ or $\Phi=\Phi_1\wedge \Phi_2$, where $\Phi_1$
and $\Phi_2$ have degree at most $d$ and level at most $l$.
If $\Phi=\neg \Phi_1$ then $\cM\vDash \Phi_1(a_1,\ldots, a_m)$ if and only
if
$\cN\vDash \Phi_1(a_1,\ldots, a_m)$, so the same holds with $\Phi$ in place
of $\Phi_1$. If $\Phi=\Phi_1\wedge \Phi_2$ then
$\cM\vDash \Phi(a_1,\ldots, a_m)$ if and only if
$\cM\vDash \Phi_i(a_1,\ldots, a_m)$, for $i=1$ and $2$, if and only if
$\cN\vDash \Phi_i(a_1,\ldots, a_m)$, for $i=1$ and $2$, if and only if
$\cN\vDash \Phi(a_1,\ldots, a_m)$. It follows that \eqref{eq:indhyp} holds
for $\Phi$ of level $l$ and any degree $d+1$; hence by induction
for all $(d,l)$.
\end{proof}

We call an expression of the form $t_1=t_2$, where $t_1$ and $t_2$
are terms, an {\em equation}. A set
 $S$  of  equations such that
every element of $S$ has variables
among $x_1,\ldots, x_m$ is called a {\em system of equations}
in $m$ variables. Let   $S$ be a system of equations in $m$ variables
and let
$\cM$ be a model of $\cL$.
We say that $(a_1,\ldots, a_m)\in M^m$ is
a {\em solution} of $S$ in $\cM$ if $\cM\vDash s(a_1,\ldots, a_m)$, for all
$s\in S$. The {\em variety} defined by $S$ over $\cM$ is the set
$V_\cM(S)=\{(a_1,\ldots, a_m)\in M^m: (a_1,\ldots,a_m)
\textrm{ is a solution of }
S\}.$ We say that a model $\cM$ of $\cL$ is {\em equationally Noetherian}
if every system $S$ of equations contains a finite subset $S_0$
such that $V_\cM(S_0)=V_\cM(S)$. As in \cite{AG1} we have the following lemma.
\begin{lem}\label{lem:eqnoeth}
Let $\cM$ and $\cN$ be $\cL$-structures.
Then,
\begin{enumerate}
\item if $\cM$ is equationally Noetherian and {\rm $\The(\cN)\subseteq\The(\cM)$}, then $\cN$ is equationally Noetherian;
\item if $\cM$ and $\cN$ are universally equivalent, $\cM$ is equationally Noetherian if and only if $\cN$ is equationally Noetherian.
\end{enumerate}
\end{lem}
\begin{proof}
Suppose that $\cM$ is equationally Noetherian and that $S$ is a system of
equations in $m$ variables.
Choose a subset $S_0\subseteq S$ such that $V_\cM(S)=V_\cM(S_0)$.
Let $S_0=\{s_1,\ldots, s_r\}$ and for each $s\in S$ let $\Phi_s$
be the sentence
$\forall x_1,\ldots ,x_m(s_1\wedge \cdots \wedge s_r\rightarrow s)$.
Since $V_\cM(S_0)=V_\cM(S)$ we have $\cM\vDash \Phi_s$ and therefore, since under the assumptions of any of the two statements above $\The(\cN)\subseteq\The(\cM)$, we have $\cN\vDash \Phi_s$, for all $s\in S$.
As $S_0\subseteq S$ it follows that $V_\cN(S)\subseteq V_\cN(S_0)$.
If $(b_1,\ldots, b_m)\in V_\cN(S_0)$ then, as $\cN\vDash \Phi_s$,
we have $(b_1,\ldots, b_m)\in V_\cN(S)$, so $V_\cN(S_0)=V_\cN(S)$.
\end{proof}
\section{Groups and Pregroups}

The language of pregroups $\Lp$ has signature $(C,F,R)$ where $C$ consists
of a single element $1$,
$F$ consists of 
a unary function symbol $^{-1}$
 and $R$ consists
of a binary relation $D$ and a
ternary relation $M$. (The usual definition of
a pregroup involves a product function defined on a subset
$D\subset P\times P$. Our description of language does not allow $F$ to
contain
partially defined functions, so we use the relation
$M$ instead of this product. We keep the relation $D$ for compatibility
with the usual definition.) A pregroup is a model $P$ of $\Lp$ satisfying
the following axioms.
\be[(i)]
\item\label{it:pre1} $\forall x,y,z((x,y,z)\in M\rightarrow (x,y)\in D)$.
\item\label{it:pre2} $\forall x,y((x,y)\in D\rightarrow \exists z((x,y,z)\in M))$.
\item\label{it:pre3} $\forall w,x,y,z((w,x,y)\in M\wedge (w,x,z)\in M\rightarrow y=z)$.
\item\label{it:pre4} $\forall x((x,1,x)\in M\wedge (1,x,x)\in M)$.
\item\label{it:pre5} $\forall x((x,x^{-1},1)\in M \wedge (x^{-1},x,1)\in M)$.
\item\label{it:pre6} $\forall x,y,z((x,y,z)\in M \rightarrow (y^{-1},x^{-1},z^{-1})\in M).$
\item\label{it:pre7} $\forall a,b,c,r,s,x((a,b,r)\in M\wedge (b,c,s)\in M \rightarrow
(
(a,s,x)\in M \leftrightarrow (r,c,x)\in M
) )$.
\item\label{it:pre8} $\forall a,b,c,d,x,y,z
(
(a,b,x)\in M\wedge (b,c,y)\in M\wedge (c,d,z)\in M
\rightarrow
\exists r,s
(
(a,y,r)\in M\vee (y,d,s)\in M
)
)$.
\ee
A pregroup homomorphism is a morphism of $\Lp$-structures
and a subpregroup is an $\Lp$-substructure of an $\Lp$-structure.
Thus $K$ is  a subpregroup of $P$ if and only if
$K$ is a pregroup, $K\subseteq P$, $1_K=1_P$, $D_K=D_P\cap (K\times K)$
 and $M_K=M_P\cap(K\times K\times K)$  (from which it follows that
the operation of inversion in $P$ extends that in $K$).

We wish, as in \cite{AG1} for the group case, to consider pregroups which
contain designated copies of some fixed pregroups (or some of their subsets).
To this end we make the following definition. 
\begin{defn}
Let $M$ be an $\cL$-structure and $N$ a subset of $M$. The {\em diagram}
of $N$ is the set of all closed atomic formulas, and their negations, which
hold in $N$.
\end{defn} 
Now let $\cS^\prime$ be a fixed multiset of  pregroups and,
for each $L\in \cS^\prime$, let  $K_L$ be a subset  of $L$ containing $1_L$.
Let $\cS$ be the set $\{K_L|L\in \cS^\prime\}$.
We define the language of $\cS$-pregroups
$\Gp$ to be the extension of $\Lp$ with signature identical to
$\Lp$ except that $C=\cup_{K\in \cS}\{d^K_k|k\in K\}$  
and $R$ contains a unary relation $\d_\cS$. A
$K${\em -pregroup} is a model $P$ of $\Gp$  satisfying the axioms
for a pregroup all the formulas of the diagram of $K$, for all $K\in \cS$,
 and 
 the additional axioms
\be[(i)]
\setcounter{enumi}{8}
\item\label{it:Gpre6}
$d^K_k\in \d_\cS$, for all $k\in K$, for all $K\in \cS$, and 
\item\label{it:Gpre7}
$\forall x(x\notin \d_\cS \rightarrow x\neq d_k)$, for all $k\in K$, 
for all $K\in \cS$.
\ee
(There is one axiom of type \eqref{it:Gpre6} and one of type
\eqref{it:Gpre7} for each $k\in K$ and $K\in \cS$
.)
A $\cS$-pregroup homomorphism is a morphism of $\Gp$-structures
and a $\cS$-subpregroup is an $\Gp$-substructure of an $\Gp$-structure.
A $\cS$-pregroup is finitely generated if it is finitely generated
as a model of   $\Gp$. If $\cS$ consists of a single element $K$ we
call an $\cS$-pregroup a $K$-pregroup and write $\Kp$ instead of $\Gp$.

\begin{lem}\label{lem:abc}
Let $P$ be a pregroup and $a,b,c\in P$. If $(a,b,c)\in M$ then
$(c,b^{-1},a)$ and $(c^{-1},a,b^{-1})\in M$.
\end{lem}
\begin{proof}
We have $(a,b,c)$ and $(b,b^{-1},1)\in M$ and, as also $(a,1,a)\in M$, axiom
\eqref{it:pre7} implies $(c,b^{-1},a)\in M$. Repeating this
argument starting with $(c^{-1},c,1)$, $(c,b^{-1},a)$ and $(1,b^{-1},b)$
 we see that $(c^{-1},a,b^{-1})\in M$.
\end{proof}

Let $\cS^\prime$ be a fixed multiset of groups and, 
for each $G\in \cS^\prime$, let 
$K_G$ be a subset of $G$ containing $1_G$. 
Let $\cS$ be the set $\{K_G|G\in \cS^\prime\}$.
 The language of $\cS$-groups is defined to be
 the language
$\cL_\cS$ with signature $(C,F,R)$, where 
$C=\cup_{K\in \cS}\{d^K_k|k\in K\}$, $F$ consists of a binary
function symbol $\cdot$ and a unary function symbol $^{-1}$ and $R$ consists
of a unary relation symbol $\d_\cS$. Then an $\cS$-group $H$ 
is a model of $\cL_\cS$
satisfying the usual group axioms with respect to $\cdot$ as multiplication 
and 
$^{-1}$ as inverse in $H$,
as well as all the formulas of the diagram of $K$, for all $K\in \cS$, 
 and the additional axioms
\be[(a)]
\item\label{it:G1}
 $d^K_k\in \d_\cS$, for all $k\in K$, $K\in \cS$, and 
\item $\forall x(x\notin \d_\S \implies x\neq d^K_k)$, for all $k\in K$,
 $K\in \cS$.
\ee
 The class of all $\cS$-groups together with the naturally defined
$\cS$-morphisms forms a category. 

If $\cS$ consists of a single element
$K$ then we refer to $K$-groups instead of  $\cS$-groups and write
$\cL_K$ instead of $\cL_\cS$. In this case, 
 if $K=G$
we recover the definition of $G$-group in \cite{AG1}.
 Further, if 
 $G=K=1$ then we drop the predicate $\d_\cS$ from the language and 
we have the standard language $\cL$ of groups. 
Note that, if $G$ is a group, a $G$-group $H$ is equationally
Noetherian in the sense defined
in the previous section if and only if it is $G$-equationally Noetherian
in the sense of \cite{AG1}.

 Notions of universal equivalence, elementary
equivalence and equivalence of finite subsets for $\cS$-groups are defined
with respect to the language $\cL_\cS$; as are substructures and extensions
of $\cS$-groups. A $\cS$-group $H$ is locally $\cS$-discriminated by a $\cS$-group
$N$ if, given a finite subset $F=\{h_1,\ldots, h_k\}$ of $H$ there is a
$\cS$-homomorphism (i.e. $\cL_\cS$-morphism) from $H$ to $N$ which is injective
on $F$. A $\cS$-group $H$ is said to be {\em finitely generated} 
if there exists
a finite subset $F$ of $H$ such that $H$ is generated by 
$F\cup \cup_{K\in \cS}K$. (Thus
a finitely generated $\cS$-group is a finitely generated $\cL_\cS$-model.)
If $P$ is any property
then a $\cS$-group $H$ is said to be locally $P$ if every 
non-trivial finitely generated
$\cS$-subgroup of $H$ has property $P$. The following theorem is proved
in  \cite{AG1}.
\begin{thm}[\cite{AG1}]\label{thm:equivgp}
Let $G$ be a group and $H$ and $K$ be $G$-groups one of which is
$G$-equationally Noetherian.
Then $H$ is locally $G$-discriminated by $K$ and $K$ is locally
$G$-discriminated by $G$
if and only if $K$ and $H$ are universally equivalent
(with respect to $\cL_G$).
\end{thm}

If $a,b$ are elements of pregroup $P$ and $(a,b)\in D_P$
we write $ab$ for the unique element $c$ such that $(a,b,c)\in M$.
Following Stallings \cite{S} we define a {\em word} of {\em length}
$k$ over a pregroup
$P$ to be a finite sequence $(c_1,\ldots, c_k)$  of elements
of $P$.
If $(c_i,c_{i+1})\in D$ then $c_ic_{i+1}\in P$ and the word
$(c_1,\ldots ,c_{i-1},c_ic_{i+1},c_{i+1},\ldots, c_k)$ is said to
be a {\em reduction} of $(c_1,\ldots, c_k)$.
The  word  $(c_1,\ldots, c_k)$ is said to be {\em reduced} if
$(c_i,c_{i+1})\notin D$, for $i=1,\ldots, k-1$.

Let $\mathbf c=(c_1,\ldots ,c_k)$ and $\mathbf a=(a_1,\ldots, a_{k-1})$
be  words such that
$(c_1,a_1)\in D$, $(a_{i-1}^{-1},c_{i})$ and $(a_{i-1}^{-1}c_{i},a_i)$
are in $D$, for $i=1,\ldots ,k-1$, and $(a_{k-1},c_k)\in D$.
 Then the {\em interleaving} $\mathbf c * \mathbf a $ of $\mathbf c$
and  $\mathbf a$ is
the word $(d_1,\ldots,d_k)$ given by
 $d_1=c_1a_1$, $d_i=a_{i-1}^{-1}c_{i}a_i$, for $i=1,\ldots ,k-1$,
and $d_k=a_{k-1}c_k$.
We define a relation $\approx$ on the set of words by
$\mathbf c\approx \mathbf d$ if and only if $\mathbf d=\mathbf c * \mathbf a$,
for some word $\mathbf a$.
As shown in \cite{S} if $\mathbf c$ is
reduced then so is $\mathbf c *\mathbf a$ and  the relation
$\approx$  is an equivalence relation on the set
of reduced words over $P$.
The {\em universal group} $U(P)$
of the pregroup $P$ is the set of equivalence classes of reduced words:
the group operation being concatenation of words followed by reduction
 to a reduced word.
As $P$ embeds in $U(P)$ then, if $P$ is a $K$-pregroup it follows that
$U(P)$ is a $K$-group.
A group $G$ may be regarded as a pregroup:
with $D=G\times G$ and $M$ the multiplication table of $G$.
It is shown in \cite{S} that $U(P)$ is universal in the
sense that, given a group $H$ and a pregroup morphism $\theta$ from $P$ to $H$,
there is a unique extension of $\theta$ to a group homomorphism from $U(P)$ to
$H$.

\begin{lem}\label{lem:preeq}
Let $P$ be a pregroup and let $(c_1,\ldots ,c_m)$ and
$(d_1,\ldots, d_n)$ be words. Then
$(c_1,\ldots ,c_m)\approx (d_1,\ldots, d_n)$ if and only if $m=n$ and
\[(d_{r-1}^{-1}\cdots d_1^{-1}c_1\cdots c_{r-1},c_r)\in D_P
\textrm{ and }
(d_r^{-1}, d_{r-1}^{-1}\cdots d_1^{-1}c_1\cdots c_{r})\in D_P,
\]
$r=1,\ldots m$, and $d_{m}^{-1}\cdots d_1^{-1}c_1\cdots c_{m}=1$.
\end{lem}
\begin{proof}
Write $D=D_P$.
Suppose first that $(c_1,\ldots ,c_m)\approx (d_1,\ldots, d_n)$.
Then by definition $m=n$ and
there exists an interleaving $(c_1,\ldots ,c_m)*(a_1,\ldots,a_{m-1})
=(d_1,\ldots, d_m)$, for some $a_i\in P$.
Then (by definition again) with $a_0=a_m=1$ we have
$(a_{i-1},c_i)$ and $(a_{i-1},c_ia_i)$ in $D$ and
$d_i=a^{i-1}c_ia_i$. Thus $(c_1,a_1)\in D$ and $d_1=c_1a_1$. Lemma
\ref{lem:abc} implies that $(d_1^{-1},c_1)\in D$ and $d_1^{-1}c_1=a_1^{-1}$.

Assume inductively that
\[
(d_{r-1}^{-1}\cdots d_1^{-1}c_1\cdots c_{r-1},c_r)\in D\textrm{ and
 }(d_r^{-1}, d_{r-1}^{-1}\cdots d_1^{-1}c_1\cdots c_{r})\in D
\]
and  $d_r^{-1}\cdots d_1^{-1}c_1\cdots c_{r}=a_r^{-1}$.
As $(a_r^{-1},c_{r+1})$ and $(a_r^{-1}c_{r+1}, a_{r+1})\in D$ and
$(a_r^{-1}c_{r+1})a_{r+1}=d_{r+1}$, Lemma \ref{lem:abc} implies
$d_{r+1}^{-1}(a_r^{-1}c_{r+1})=a_{r+1}^{-1}$. Combined with the
inductive hypothesis this shows that the $(r+1)$st version of this hypothesis
also holds. Hence the statement of the inductive hypothesis holds
for $r=1,\ldots ,m$. As $a_m=1$ we obtain, from the $k$th version
 $d_{m}^{-1}\cdots d_1^{-1}c_1\cdots c_{m}=1$, as required.

Conversely, suppose the conditions given in the lemma hold. Then
$(d_1^{-1},c_1)\in D$ and so we may define $a_1^{-1}=d_1^{-1}c_1$. Two
applications of Lemma \ref{lem:abc} show that $(c_1,a_1)\in D$ and
$c_1a_1=d_1$. Define $a_0=1$ and suppose that $a_1,\ldots ,a_r$ have
been defined such that $(a_{i-1}^{-1},c_i), (a_{i-1}^{-1}c_i,a_i)\in D$
$a_i^{-1}=d_i^{-1}\cdots d_1^{-1}c_1\cdots d_i$ and $d_i=a_{i-1}^{-1}d_ia_i$,
$i=1,\ldots r$. Then $(d_r^{-1}\cdots d_1^{-1}c_1\cdots c_r, c_{r+1})$ and
$(d_{r+1}^{-1}, d_r^{-1}\cdots d_1^{-1}c_1\cdots c_{r+1})\in D$ and we
may set $a_{r+1}^{-1}= d_{r+1}^{-1}\cdots d_1^{-1}c_1\cdots c_{r+1}
=d_{r+1}^{-1}(a_r^{-1}c_{r+1})$. Two applications of Lemma \ref{lem:abc}
give $a_r^{-1}c_{r+1}a_{r+1}=d_{r+1}$. Finally we obtain
$a_m^{-1}=d_{m}^{-1}\cdots d_1^{-1}c_1\cdots c_{m}=1$ so
$(c_1,\ldots ,c_m)*(a_1,\ldots, a_{m-1})$ is defined and equal to
$(d_1,\ldots, d_m)$ as required.
\end{proof}
\begin{cor}\label{cor:subpre}
If $Q$ is a subpregroup of a pregroup $P$ then $U(Q)$ is
a subgroup of $U(P)$. In particular, if $P$ is an $\cS$-pregroup
 then 
$U(P)$ is an $\cS$-group.
\end{cor}
\begin{proof}
To prove the first statement we 
need to show that if $\mbf a$ and $\mbf b$ are words over $Q$ then
$\mbf a\approx \mbf b$ in $Q$ if and only if $\mbf a\approx\mbf b$
in $P$. Suppose that $\mbf a\approx  \mbf b$ in $P$. Then
using Lemma \ref{lem:preeq} and the definition of $\Lp$-substructure
we have $\mbf a\approx  \mbf b$ in $Q$. As the opposite implication
is immediate this proves the first part of the corollary. For the second
statement suppose that $K\in \cS$ and that $K$ is a subset
of a pregroup $L$, as in the definition above. As $K\subseteq P$ we 
may assume that $L\subseteq P$ and so $K\subseteq U(L)\subseteq U(P)$.
\end{proof}
\begin{thm}\label{thm:ugequiv}
Let
$P_1$ and $P_2$ be $\cS$-pregroups.
If $P_1\equiv_\exists
P_2$ with respect to $\Lp_\cS$ then $U(P_1)\equiv_\exists U(P_2)$
with respect to $\cL_\cS$.
\end{thm}
\begin{proof}
Let $U_i=U(P_i)$ and $D_i=D_{P_i}$,  for $i=1,2$.
 We shall show that $\cF(U_1)\equiv \cF(U_2)$
and the theorem will then follow from Proposition \ref{prop:fe}.

Let $F=\{\tilde u_1,\ldots ,\tilde u_m\}$ be a finite subset of $U_1$.
For each $i$ choose a representative $u_i$ of $\tilde u_i$
and write it as a reduced word $u_i=(c_{i1},\ldots ,c_{im_i})$
over $P_1$. Let $S_0=\cup_{i=1}^m\cup_{j=1}^{m_i}
\{c_{ij}\}$ and for all $r\ge 0$ let $S_{r+1}=S_r\cup
\{ab:a,b\in S_i\textrm{ and } (a,b)\in D_1\}$.
Let $J=\max\{m_i:i=1,\ldots ,m\}$ and define $S=S_{2J}$. As $P_1\equiv_\exists
P_2$ there is, using Proposition \ref{prop:fe},
an $\Lp_\cS$-isomorphism $\phi$ from $S$ to a subset $T$ of $P_2$.
Note that setting $T_0=\phi(S_0)$ we may define $T_r$ as we  have
defined $S_r$, with $T_0$ in place of $S_0$ and $P_2$ in
place of $P_1$. Then, by definition of isomorphism
and by construction of $S$ it follows that $\phi(S_r)=T_r$,
for $r=0,\ldots, 2J$,
so $T=T_{2J}$.
Let $\mathbf c=(c_1,\ldots, c_k)$ be a word over $S$
(i.e. $c_i\in S$, for all $i$)
and let $\phi(c_i)=d_i$. Then $\mathbf d=(d_1,\ldots, d_k)$ is a word over $T$ and
we define a map $\theta$ from words over $S$ to words over $T$ by
setting $\theta(\mathbf c)=\mathbf d$. In this case, for all such $\mbf c$,
we have $(c_i,c_{i+1})\in D_1$ if and only if $(d_i,d_{i+1})\in D_2$, so
$\theta$ maps reduced words to reduced words.

Now let $\mbf p_1=(p_{11},\ldots ,p_{1m_1})$ and
$\mbf p_2=(p_{21},\ldots, p_{2m_2})$ be words over $S_0$, with
$m_i\le J$. Let $\phi(p_{ij})=q_{ij}$, and let
$\theta(\mbf p_i)=\mbf q_i=(q_{i1},\ldots ,q_{im_i})$, $i=1,2$.
From Lemma \ref{lem:preeq} we have $\mbf p_1\approx\mbf p_2$ if and
only if $m_1=m_2=k$, $(p_{2,r-1}^{-1}\cdots p_{2,1}^{-1}
p_{1,1}\cdots p_{1,r-1},
p_{1,r})$ and $(p_{2,r}^{-1},p_{2,r-1}^{-1}\cdots p_{2,1}^{-1}
p_{1,1}\cdots p_{1,r-1}p_{1,r})$ belong to $D_1$, for $r=1,\ldots ,k$, and
$ p_{2,k}^{-1}\cdots p_{2,1}^{-1}
p_{1,1}\cdots p_{1,k}=1$. Since all the elements of $P_1$ involved in
these conditions belong to $S$, the conditions hold if and only if
they hold on replacing $p_{ij}$ with $q_{ij}$. Hence
$\mbf p_1\approx\mbf p_2$ if and only if $\mbf q_1\approx\mbf q_2$.
Therefore $\theta$ induces a map $\tilde \theta$ from equivalence classes of reduced
words over $S_0$, of length at most $J$, to equivalence classes of
reduced words over $T_0$.

Let $\tilde S$ and $\tilde T$ denote the sets of
equivalence
classes of reduced words of length at most $J$, over $S_0$
and $T_0$ respectively.
To see that the map that $\tilde \theta$  is an $\cL_\cS$-morphism
from $\tilde S$ to $\tilde T$
consider
a word (not necessarily reduced) $\mathbf p=(p_1,\ldots, p_k)$ over $S_0$
of length $k\le J$.
Let $q_i=\phi(p_i)$ and let $\theta(\mbf p)=\mbf q=(q_1,\ldots, q_k)$.
We claim that for $r$ with $0\le r\le k-1$ there is a sequence
of $r$ reductions which we may apply
to $\mbf p$, resulting in a word $\mbf p_r$, if and only if
there is a corresponding sequence of $r$ reductions which
we may apply  to $\mbf q$ resulting in a word $\mbf q_r$ such that
$\theta(\mbf p_r)=\mbf q_r$. Moreover $\mbf p_r\in S_r$ and
$\mbf q_r\in T_r$. This holds trivially for $r=0$.
Suppose that it holds for $0,\ldots ,r$, for some $0\le r\le k-2$.
Let $\mbf p_r=(p_{r,1},\ldots,p_{r,s})$ and
$\mbf q_r=((q_{r,1},\ldots,q_{r,s})$, with $\mbf p_r\in S_r$ and
$\mbf q_r\in T_r$ and $\mbf q_r=\theta(\mbf p_r)$. We may apply a reduction to $\mbf p_r$ if and
only if $(p_{r,i},p_{r,i+1})\in D_1$, for some $i$, in which case we
may define $\mbf p_{r+1}=(p_{r,1},\ldots,p_{r,i}p_{r,i+1},\ldots p_{r,s})$
and then $\mbf p_{r+1}\in S_{r+1}$. Since $\phi$ is an $\Lp_\cS$-isomorphism
this occurs if and only if $(q_{r,i},q_{r,i+1})\in D_2$, in which case we
may define $\mbf q_{r+1}=(q_{r,1},\ldots,q_{r,i}q_{r,i+1},\ldots q_{r,s})$
and then $\mbf q_{r+1}\in T_{r+1}$. Since $\theta(\mbf p_r)=\mbf q_r$
 it follows that $\theta(\mbf p_{r+1})=\mbf q_{r+1}$ and so
the claim holds for all $r$. Now let $\tilde p_1$ and $\tilde p_2$ be
elements of $\tilde S$ and let $\mbf p_1$ and $\mbf p_2$ be reduced
words, of length at most $J$, over $S_0$ representing $\tilde p_1$ and
$\tilde p_2$, respectively. Suppose that $\mbf p$ is a reduced word
(over $S$)
obtained from the concatenation $\mbf p_1 \mbf p_2$ by a sequence of
reductions. Then, in $U_1$, we have $\tilde p_1\tilde p_2=\tilde p$, where
$\tilde p$ is the equivalence class of $\mbf p$.
Let $\mbf q_i=\theta (\mbf p_i)$ and $\mbf q=\theta (\mbf p)$. Then,
from the above, the concatenation $\mbf q_1\mbf q_2$ reduces to $\mbf q$,
which is a reduced word over $T$.
Hence, in $U_2$, $\tilde q_1\tilde q_2=\tilde q$, where $\tilde q$ is the
equivalence class of $\mbf q$. Now, in the case where $\mbf p$ is a word
over $S_0$ we have $\tilde\theta(\tilde p_1)\tilde\theta(\tilde p_2)
=\tilde q_1\tilde q_2=\tilde q=\tilde\theta(\tilde p)
=\tilde\theta(\tilde p_1\tilde p_2)$, showing that $\tilde \theta$ is
an $\cL_\cS$-morphism. Using the result of the first half of this paragraph
and the fact that $\tilde \theta$ is bijective we can show that
$\tilde \theta^{-1}$ is  also an $\cL_\cS$-morphism. In particular
$\tilde\theta$ restricted to $F$ is an $\cL_\cS$ isomorphism onto its image.
Therefore $\cF(U_1)=\cF(U_2)$, as required.
\end{proof}

\section{Applications}

In this section we apply Theorem \ref{thm:ugequiv} to prove that the universal equivalence of pregroups translates nicely into universal equivalence of free constructions.

\subsection{Free products}
To simplify notation we assume from the outset that we have 
two groups $A$ and $B$ whose intersection is the identity element. 
In this case let $P=A\cup B$ 
and set $D=(A\times A)\cup (B\times B)$. Then $P$ is a pregroup and 
$U(P)=A\ast B$.

\begin{prop}\label{prop:frpr}
  Let $A_1, B_1, A_2$ and $B_2$ be groups such that 
$A_1\cap B_1=A_2\cap B_2=1$.  If $\cF(A_1)\equiv \cF(A_2)$ and $\cF(B_1)\equiv \cF(B_2)$ then $A_1\ast B_1$ is existentially equivalent to $A_2\ast B_2$.
\end{prop}
\begin{proof}
 Let $P_1=A_1\cup B_1$ and $P_2=A_2\cup B_2$ be two pregroups as above.
  Let $S$ be a finite subset of $P_1$ in the language $\Lp$. Then $S=(S\cap A_1)\cup (S\cap B_1)=S_{A_1}\cup S_{B_1}$. Let $S_{A_2}'$ and $S_{B_2}'$ be two finite subsets of $A_2$ and $B_2$ in the language of groups $\cL$, isomorphic to $S_{A_1}$ and $S_{B_1}$, respectively. Then $S'=S_{A_2}'\cup S_{B_2}'$ is a subset of $P_2$ isomorphic to $S$ in the language $\Lp$. By Proposition \ref{prop:fe}, $P_1\equiv_\exists P_2$ in the language $\Lp$, and by Theorem \ref{thm:ugequiv} $A_1*B_1\equiv_\exists A_2*B_2$ in the language $\cL$.
\end{proof}

\subsection{Free Products with Amalgamation}
Again it simplifies notation to assume that $A$ and $B$ are $C$-groups
which intersect in the designated copy of the subgroup $C$, where 
$C\neq 1$. In this case 
let $P=A\cup B$
and set $D=(A\times A)\cup (B\times B)$. Then $P$ is a $C$-pregroup and
$U(P)=A*_C B$. 

\begin{prop} \label{prop:amalg}
 Let $A_1, B_1, A_2$ and $B_2$ be $C$-groups such that 
$A_1\cap B_1=A_2\cap B_2=C$.
 If $\cF_{\cL_C}(A_1)\equiv \cF_{\cL_C}(A_2)$ and 
$\cF_{\cL_C}(B_1)\equiv \cF_{\cL_C}(B_2)$ then the group $A_1\ast_{C_1} B_1$ is existentially equivalent to $A_2\ast_{C_2}B_2$ in the language $\cL_C$ and,
a fortiori, in the language $\cL$.
\end{prop}
\begin{proof}
Let $P_1=A_1\cup B_1$ and $P_2=A_2\cup B_2$, 
be two $C$-pregroups as above.
Let $S$ be a finite subset of $P_1$ in the language $\Lp_C$.  Let
$S_{A_1}=S\cap A_1$ and $S_{B_1}=S\cap B_1$ so 
 $S=S_{A_1}\cup S_{B_1}$. Let $S_{A_2}'$ and $S_{B_2}'$ be two finite subsets of $A_2$ and $B_2$ in the language of $C$-groups $\cL_C$, isomorphic to $S_{A_1}$ and $S_{B_1}$, respectively. Then $S'=S_{A_2}'\cup S_{B_2}'$ is a subset of $P_2$ isomorphic to $S$ in the language $\Lp_C$. By Proposition \ref{prop:fe}, $P_1\equiv_\exists P_2$ in the language $\Lp_C$, and by Theorem \ref{thm:ugequiv} $A_1*_CB_1\equiv_\exists A_2*_C  B_2$ in the language $\cL_C$.
\end{proof}

\subsection{HNN-Extensions}

Given a group $G$ and an isomorphism $\theta:C_1\rightarrow C_2$, where
$C_1$ and $C_2$ are subgroups of $G$, let $t$ be a symbol not in $G$ and
\begin{equation} \label{eq:HNN}
P_0=G\cup t^{-1}G\cup Gt \cup t^{-1}Gt.
\end{equation}
Let $P$ be the set of equivalence classes of the equivalence relation generated by $t^{-1}ht=\theta(h)$, for all $h\in C_1$.
Set
\[
D=\bigcup_{
\e_0,\e_1,\e_2=0,1
}t^{-\e_0}Gt^{\e_1}\times t^{-\e_1}Gt^{\e_2}\subseteq P\times P.
\]
(The equivalence 
relation means that $(p,\theta(c))$ and $(\theta(c),p)$ belong to 
$D$ for all $c\in C_1$ and $p\in P$.) Then $P$ and $D$ constitute
a pregroup and it can be verified that $C_1$ and $C_2$ embed in $P$.
Hence $P$ is an $\cS$-pregroup, where $\cS=\{C_1, C_2\}$. 
Moreover $U(P)$ is 
the HNN-extension
$\la G, t| t^{-1}ct=\theta(c), c\in C_1\ra$, which  is an $\cS$-group
(with constants $C_1$ and $\theta(C_1)=C_2$). 

\begin{prop} \label{prop:HNN}
Let $A_1$ and $A_2$ be $\cS$-groups, where $\cS=\{C_1,C_2\}$ and 
$\theta:C_1\rightarrow C_2$ is an isomorphism.
 If $\cF_{\cL_{\cS}}(A_1)\equiv \cF_{\cL_{\cS}}(A_2)$ 
then the group $G_1=\langle A_1, t\mid t^{-1}ct=\theta(c), c\in C_1\rangle$ 
is existentially equivalent, in the language $\cL_\cS$ and in 
the language $\cL$, to the group 
$G_2=\langle A_2, t\mid t^{-1}ct=c, c\in C_1\rangle$.
\end{prop}
\begin{proof}
Let $P_1$ and $P_2$ be the  two $\cS$-pregroups corresponding
to $A_1$ and $A_2$, respectively, as defined above, and let
$P_{1,0}$ and $P_{2,0}$ be the underlying sets, as in \eqref{eq:HNN}.
Let $S$ be a finite subset of $P_1$ in the language $\Lp_{\cS}$. 
Let $\hat S\subseteq P_{1,0}$ be the union of all the equivalence classes
of elements of $S$. 
Then $\hat S$ is a disjoint union of $4$ sets, $S_1=\hat S\cap A_1$, 
$S_2=\hat S\cap t^{-1}A_1$,  $S_3=\hat S \cap A_1 t$ and 
$S_4=\hat S\cap t^{-1}A_1t$.
To obtain corresponding sets in $A_1$ define 
$T_1= S_1$, $T_2=t S_2$, $T_3=S_3 t^{-1}$ and $T_4=t S_4 t^{-1}$. 
By hypothesis there exist subsets $T_i^\prime \subseteq A_2$, such that
$T_i\cong_{\cL_\cS} T_i^\prime$, for $i=1,\ldots , 4$. Set
$S_1^\prime=T_1^\prime$, $S_2^\prime= t^{-1} T_2^\prime$, 
$S_3^\prime= T_3^\prime t$ and  
$S_4^\prime= t^{-1}T_4^\prime t$. Define 
$\hat S^\prime=
S_1^\prime\cup S_2^\prime\cup S_3^\prime\cup S_4^\prime$.
The $\cL_\cS$-isomorphisms between the 
$T_i$'s and the $T_i^\prime$'s induce a bijection from 
$\hat S$ to $\hat S^\prime$ and by construction this isomorphism
factors through the equivalence relations on $P_{1,0}$ and $P_{2,0}$
to give an $\Lp_\cS$-isomorphism between $S$ and the quotient $S^\prime$ of
$\hat S^\prime$ in $P_2$. Applying Proposition \ref{prop:fe} and Theorem \ref{thm:ugequiv}, $G_1$ is universally equivalent to $G_2$ in the language
 $\cL_{\cS}$ (and consequently in the language $\cL$). 
\end{proof}


\begin{thebibliography}{99}
\bibitem{AG1} G. Baumslag, A. Myasnikov, and V.N. Remeslennikov, {\it Algebraic
geometry over groups. I. Algebraic sets and ideal theory,}  J. Algebra {\bf 219} (1999), 16-79.
\bibitem{marker} D. Marker {\it Model Theory: An Introduction.} Springer, 2002.
\bibitem{S} J.~Stallings, {\em Group Theory and Three-Dimensional Manifolds.}
Yale Mathematical Monographs {\bf 4}, Yale Univ. Press, 1971.
\end{thebibliography}
\end{document}